\documentclass{article}
\usepackage{amsmath,amsthm,amsfonts,amscd,amssymb,eucal,latexsym,hyperref} 
\usepackage[all]{xy}
\begin{document}

\title{Growth rates of dimensional invariants of compact quantum groups and a  theorem of H\o egh-Krohn, Landstad and St\o rmer}
  
\author{Claudia Pinzari}
\date{}
\maketitle

   \centerline{\it Dedicated to the memory of Claudio D'Antoni}                
\bigskip

\begin{abstract}
We give local upper and lower bounds for the eigenvalues of the modular operator associated 
to an ergodic action of a compact quantum group on a unital $C^*$--algebra. They involve the modular theory of the quantum group and
the growth rate of quantum dimensions of its representations and they become sharp 
if other integral invariants grow subexponentially. For compact groups, this reduces to the finiteness theorem  of  H\o egh-Krohn, Landstad and St\o rmer. Consequently,  compact quantum groups of Kac type admitting an ergodic action with a non-tracial invariant state must have representations whose dimensions grow exponentially. 
In particular, $S_{-1}U(d)$ acts ergodically only on tracial $C^*$--algebras.
For   quantum groups   with non-involutive coinverse, we derive a lower bound  for the
parameters $0<\lambda<1$ of factors of  type III${}_\lambda$ that can possibly arise from the  GNS representation of the invariant state  of an ergodic action with a factorial centralizer.
\end{abstract}

\section{Introduction}
In the early 80's H\o egh-Krohn, Landstad and St\o rmer proved  that the multiplicity
of an irreducible representation of a compact group acting ergodically on a unital $C^*$--algebra   is bounded above by its dimension and moreover the unique invariant state  is a trace  \cite{HLS}.  This   result is often used as a finiteness criterion in operator algebras. Moreover,    Wassermann, starting from this, showed the negative result that   ${\rm SU}(2)$ acts ergodically only on type I von Neumann algebras \cite{Wassermann3}.  

 If we consider compact quantum groups instead of compact groups,  ergodic theory on operator algebras becomes much richer. For example, as is
 well known, finiteness   fails, as compact quantum groups may have non-involutive coinverse 
and, in this case, the Haar  state  has a non-trivial modular theory
 \cite{WoronowiczCMP}.
Boca generalized some of these results to ergodic actions of compact  quantum groups.  He proved   that the quantum dimension is
 an  upper bound  for the multiplicity and that the invariant state     satisfies the KMS  property \cite{Boca}. 
Wang has found many examples of ergodic actions of the free unitary quantum groups on factors of type II and III  \cite{Wang}. 
In particular, his examples show that finiteness fails already in the case of compact quantum groups with involutive coinverse, often called of Kac type,    even though
 these have trivial modular theory. Indeed, $A_u(n)$ can act ergodically 
 on factors  on type III${}_\frac{1}{n}$. 
Bichon, De Rijdt and Vaes have constructed actions of $A_o(F)$
with multiplicities larger than the integral dimensions  \cite{BDV}. In joint papers with Roberts, we have shown   that any finite index inclusion of factors of  type II${}_1$
gives rise to ergodic actions of $A_o(F)$  \cite{PRsubfactors}. Moreover,  the classification of ergodic actions of compact quantum groups  is related to the classification of tensor $C^*$--categories with conjugates \cite{PRinduction}. 

The aim of this note is to show an analogue of the second part of the finiteness theorem of \cite{HLS} for compact quantum groups.

The examples of \cite{BDV}, as well as subsequent developments of the general theory
of ergodic actions \cite{PR}, contributed to a further understanding of the invariant state. It became clear, for example, that it  is always
almost periodic in the sense of Connes \cite{ConnesAP}, meaning that the associated modular operator is diagonal, and that its eigenvalues  exhibit an explicit, separated, dependence   on both the  modular theory of the quantum group and the eigenvalues of certain strictly positive matrices canonically associated
to the irreducible represenations of the quantum group belonging in the spectrum of the action. 

A consequence of this is that, in the special case where
the invariant state of an action  of a non-Kac type compact quantum group  is factorial and has factorial centralizer, then the associated factor is necessarily
of type III${}_\lambda$ with $0<\lambda\leq1$  (cf. Sect. 3 for a more precise statement). Note that  the invariant state of the mentioned examples of Wang   satisfy this factoriality condition. Another class of examples has been  recently studied by
Vaes and Vergnioux   \cite{VV}. They  considered the translation action of $A_o(F)$ over itself and they
were able to prove, among other things, that if $F$ has rank at least $3$ and satisfies a suitable  condition, then the Haar state is factorial
and moreover the factor generated by the the GNS construction is  full in the sense of Connes \cite{ConnesAP}.

For a general ergodic action, the problem of studying modularity of the invariant state  becomes that of studying the spectra of the associated  
 matrices. In retrospect, the theorem of \cite{HLS} amounts 
to showing that these  matrices are always trivial in the classical case. 

To study these spectra, we  introduce two growth rates of dimensional invariants of the compact quantum group,  that of integral dimensions, $\text{Dim}_u$, and    of quantum   dimensions, $D_u$. 
These growth rates
 distinguish between $SU(2)$, $S_qU(2)$ for  $0<q<1$, $A_o(F)$, for $\text{rank}(F)\geq3$, and  $A_u(F)$. More precisely,
in the group case, quantum dimensions are just the integral dimensions, and therefore the 
growth rate is always polynomial, as a consequence of Weyl's dimension formula, hence $D_u=\text{Dim}_u=1$ for all representations. This polynomial growth rate
played an important role  in the original proof of the finitenss theorem of \cite{HLS}.

Growth rate of quantum dimensions for $S_qU(2)$, $0<q<1$  is instead exponential. However,  integral dimensions of irreducibles are the same as in the classical case, hence they  grow polynomially.
This is   opposed to $A_o(F)$, for $\text{rank}(F)\geq3$, for which both growth rates are exponential. In the case of $A_u(F)$, both dimensional invariants of the fundamental representation  have the largest possible growth rates as all tensor powers of this representation are irreducible. 

We shall see that if the growth rate of the integral dimensions is 
subexponential, then that of quantum dimensions,  $D_u$, is explicitly determined by the eigenvalues of the  modular operator of the quantum group (Prop. 4.3). However, this fact is not generally true, and $A_o(n)$, for $n\geq3$, or
$A_u(n)$, for $n\geq2$, are the first examples.

Given an ergodic action of a compact quantum group on a unital $C^*$--algebra, our main result consists of showing that $D_u$ and $D_u^{-1}$ are upper and lower  bounds for the  spectral radii of the
 mentioned  matrices  involved in the modular operator of the ergodic action (see Theorem 4.8).  In particular, we reproduce the finiteness theorem of \cite{HLS} from the fact that $D_u=1$ in the group case. Our bounds  become equalities
 in the case of ergodic actions with high quantum multiplicities in the sense of \cite{BDV}, but  integral multiplicities with subexponential growth rate. 
 
 We shall show that 
  the invariant $D$  determines the parameter $\lambda=\frac{1}{n}$ of the Connes invariant $S$ in Wang's example with $A_u(n)$.
In other words, the rate of growth of dimensions   clarifies which compact quantum groups of Kac type can act    on infinite factors.
More precisely, one of the consequences of our result is that, among compact quantum groups of Kac type,
only those having some irreducible representation whose integral dimension grows exponentially
can possibly act on $C^*$--algebras with a non-tracial invariant state.
For example, all ergodic $C^*$--algebras for  $S_{-1}U(d)$ are tracial. 

For compact quantum groups with non-involutive coinverse, we  derive a general lower bound for the possible parameters $0<\lambda<1$ of the 
 type III${}_\lambda$ factors than can arise from the GNS representation of the invariant state, provided the centralizer algebra is  a factor. This lower bound   involves the modular theory of the quantum group, the growth rates $D_u$, and the spectrum of the
 action,  see Cor. 4.13. We thus see that many parameters $\lambda$ may be excluded
if  some spectral information is known.
For example, for $S_qU(2)$, with $0<|q|<1$ then $\lambda\geq |q|^{2r}$, where $r$ labels the first spectral irreducible representation. For $A_o(F)$, with $F$   of rank $\geq3$,  then $\lambda\geq (\frac{q}{\|F\|^2})^r$ where $q+\frac{1}{q}:=\text{Trace}(F^*F)$ and $r$ is as before. As another example,  
$\lambda\geq\frac{\min\{q_0,q_n^{-1}\}}{\text{Trace}(F^*F)}$,  for $A_u(F)$, with  $q_0$ and $q_n$ the smallest and largest eigenvalues of $F^*F$ if the fundamental representation is spectral.  

Our methods differ from the original proof of the classical finiteness theorem. They rely on the duality theorem of \cite{PR}, which allows a rather simple presentation. 

This paper is organized as follows. Sect. 2 is devoted to the preliminaries: we recall the notion of standard solutions of the conjugate equations and  the duality theorem 
for ergodic actions. In Sect. 3 we describe preliminary  results on the modular theory of the invariant state. In Sect. 4 we discuss the main result and some corollaries.

\section{Preliminaries}

\noindent{\it Standard solutions of the conjugate equations.} We shall briefly recall the main
features of tensor $C^*$--categories with conjugates   \cite{LR}. 
Let ${\cal A}$ be a tensor $C^*$--category 
 (always assumed to be strict, with irreducible tensor unit $\iota$, subobjects and direct sums).
Arrows $R\in(\iota, \overline{u}\otimes u)$ and $\overline{R}\in(\iota, u\otimes \overline{u})$ are said to  define a conjugate $\overline{u}$ of the object $u$  if  they satisfy the conjugate equations $$\overline{R}^*\otimes 1_u\circ 1_u\otimes R=1_u,\quad {R}^*\otimes 1_{\overline{u}}\circ 1_{\overline{u}}\otimes \overline{R}=1_{\overline{u}}.$$
If $R_u$, $\overline{R}_u$ and $R_v$, $\overline{R}_v$ are solutions for $u$ and $v$ then   $R_{u\otimes v}:=1_{\overline{v}}\otimes R_u\otimes 1_v\circ R_v$, $\overline{R}_{u\otimes v}:=1_u\otimes \overline{R}_v\otimes 1_{\overline{u}}\circ \overline{R}_u$ are solutions for  $u\otimes v$. Similarly, $R_{\overline{u}}:=\overline{R}_u$, $\overline{R}_{\overline{u}}:=R_u$ are solutions for $\overline{u}$.  They are called the tensor product and  conjugate  solutions respectively. If   conjugates exist for every object, then every object is the direct sum of irreducible objects. 

The category  $\text{Rep}(G)$ of (unitary, finite dimensional) representations  of a compact quantum group has conjugates and
is embedded in the category of finite dimensional 
 Hilbert spaces, hence arrows are linear maps. Solutions of the conjugate equations take the 
 form $R=\sum_i j\psi_i\otimes\psi_i$, $\overline{R}=\sum_k j^{-1}\phi_k\otimes\phi_k$, with
 $j:H_u\to H_{\overline{u}}$ a unique invertible antilinear map between the representation Hilbert spaces and $(\psi_i)$, $(\phi_k)$ orthonormal bases. 
 Conjugate and tensor product solutions correspond to $j_{\overline{u}}:=j_{u}^{-1}$, $j_{u\otimes v}=j_v\otimes j_u\theta$, with $\theta:H_{u}\otimes H_v\to H_v\otimes H_u$ the flip map.
 
 In a tensor $C^*$--category with conjugates it is  convenient to select {\it standard solutions} of the conjugate equations, meaning that $\|R_u\|=\|\overline{R}_u\|$ if $u$ is an irreducible object, while if $u\simeq \oplus_i u_i$ with $u_i$ irreducible,     $R_u:=\sum_i \overline{S}_i\otimes S_i\circ R_{u_i}$ $\overline{R}_u:=\sum_i {S}_i\otimes \overline{S}_i\circ \overline{R}_{u_i}$, where $R_{u_i}$, $\overline{R}_{u_i}$ are standard solutions and $\{S_i\in(u_i, u)\}$, $\{\overline{S}_i\in(\overline{u}_i, \overline{u})\}$ are two sets of isometries whose ranges
are pairwise orthogonal  add up to $1_u$ and $1_{\overline{u}}$, respectively. Clearly, the conjugate of a standard solution is standard. Standard solutions are unique up to {\it unitary} equivalence. Most importantly, they realize the minimal value of $\|R_u\|\|\overline{R}_u\|$ among all solutions, and they are characterized, up to scalars, by this property. This minimal value is the quantum (or intrinsic) dimension of $u$,   denoted by $d(u)$. This implies the  following fact, which will play a role.
\medskip

\noindent{\bf 2.1. Theorem.} \cite{LR} {\it 
The   tensor product of standard solutions is  standard.}\medskip
 
We shall also need the following fact. In the category $\text{Rep}(G)$, if $j_u: H_u\to H_{\overline{u}}$ defines a standard
 solution of the conjugate equations for $u$, $j_u^*j_u$ is a positive operator on $H_u$
 which does not depend on the choice of the standard solution. This operator 
correspond to  $F_u^{-1}$ of Woronowicz \cite{WoronowiczCMP}. It follows that
 the spectrum of $j_u^*j_u$ as well as its smallest and largest eigenvalues, denoted $\lambda_u$ and $\Lambda_u$ respectively, are invariantly associated to $u$ (in fact, to its equivalence class). 
 Note that $\lambda_{\overline{u}}=\Lambda_{u}^{-1}$, $\lambda_{u\otimes v}=\lambda_u\lambda_v$, and similarly for $\Lambda_{u\otimes v}$.
 \medskip

\noindent{\it Ergodic actions, spectral functor and quasitensor functors.}
To a certain extent,  we may think of the relationship between an ergodic action of a compact quantum group on a unital $C^*$-algebra, the associated spectral functor and an abstract quasitensor functor,  as analogous to that  between a Lie group, the associated
 Lie algebra and an abstract Lie algebra. The analogy is supported by the following properties of ergodic actions proved in \cite{PR}, and referred to as a duality theorem for ergodic actions. The spectral functor of an ergodic action is a quasitensor functor. For any ergodic  action, there always is a maximal ergodic action which has the same spectral functor as that of the original one. It is the completion in the maximal $C^*$--norm of the dense spectral subalgebra of the given ergodic action. The maximal ergodic action with a given spectral functor is unique and canonically associated with it.
Any abstract quasitensor functor from the representation category of $G$
to the Hilbert spaces is the spectral functor of an ergodic action. Two maximal ergodic actions of a given compact quantum group $G$ are conjugate if and only if the associated spectral functors are related by a unitary natural transformation. Hence quasitensor functors $\text{Rep}(G)\to\text{Hilb}$ classify maximal ergodic $C^*$--actions of  $G$.

In some more   detail, let $G$ be a compact quantum group 
\cite{WoronowiczLesHouches} and
$$\alpha:{\cal C}\to {\cal C}\otimes {\cal Q}$$ an ergodic action of $G$   on a unital $C^*$--algebra ${\cal C}$ (i.e. ${\cal C}^\alpha:=\{c\in{\cal C}, \alpha(c)=c\otimes I\}={\mathbb C}$). ${\cal Q}$ denotes the Hopf $C^*$--algebra of $G$.
Let $u$ be a   representation of $G$ on the Hilbert space $H_u$.  
Consider the space   of linear maps $T:H_u\to{\cal C}$ intertwining $u$ with the  action $\alpha$. The map taking $u$ to this space of intertwiners   is a functor from the representation category of $G$ to the category of vector spaces. However, it is a contravariant functor, hence   for convenience we consider  the covariant functor   obtained passing to the category of dual vector spaces. This covariant functor will be donoted by $L$ (while we used the notation $\overline{L}$ in previous papers).
As a consequence of ergodicity, the predual of $L_u$, and hence of $L_u$ itself, is a Hilbert space. We may  explicit $L_u=\{\sum_i\psi_i\otimes c_,\quad \alpha(c_i)=\sum_jc_j\otimes u_{j,i}^*, \quad \psi_i \text{ o. b.}\}$ and thus think of $L_u$ as the space of fixed points $(H_u\otimes{\cal C})^{u\otimes\alpha}$ with inner product arising
from the restriction of  the ${\cal C}$--valued inner product of the free Hilbert module
 $H_u\otimes{\cal C}$.
  If $A\in(u,v)$ is an intertwiner in $\text{Rep}(G)$, ${L}_A$ acts as $A\otimes I$ from $L_u$ to $L_v$.

A $^*$--functor $\mu:{\cal A}\to{\cal M}$  between two tensor $C^*$--categories
is called  quasitensor if  there are   isometries
$\tilde\mu_{u,v}\in(\mu_u\otimes\mu_v, \mu_{u\otimes v})$,
 such that
$${\mu}_\iota=\iota,\eqno (2.1)$$
$$\tilde\mu_{u,\iota}=\tilde\mu_{\iota, u}=1_{\mu_u},\eqno(2.2)$$
$$\tilde\mu_{u,v\otimes w}^*\circ\tilde\mu_{u\otimes
v,w}=1_{\mu_u}\otimes\tilde\mu_{v,w}
\circ{\tilde\mu_{u,v}}^*\otimes 1_{\mu_w}\eqno(2.3)$$
and natural in $u$, $v$, 
$$\mu({S\otimes T})\circ
\tilde\mu_{u,v}=\tilde\mu_{u',v'}\circ\mu(S)\otimes\mu(T),\eqno(2.4)$$
 for  objects $u$, $v$, $w$, $u'$, $v'$ of ${\cal A}$ and arrows
$S\in(u, u')$,
$T\in(v,v')$.
This definition was given in \cite{PR}.
Note that  the most relevant  axiom, $(2.3)$, implies associativity: $\tilde{\mu}_{u,v,w}:=\tilde{\mu}_{u\otimes v, w}\circ\tilde{\mu}_{u,v}\otimes 1_{\mu_w} =\tilde{\mu}_{u, v\otimes w}\circ1_{\mu_u}\otimes\tilde{\mu}_{v,w}$.
If all the isometries $\tilde\mu_{u,v}$ are unitary, $(2.3)$, is equivalent to associativity. In this case
$(\mu,\tilde{\mu})$ will be called a {\it relaxed tensor} functor.  Examples arise from pairs
of  non-isomorphic compact quantum groups 
with tensor equivalent representation categories. A well known class of examples is 
$S_\mu U(2)$ and $A_o(F)$ for suitable conditions on $F$ (cf. Example 4.6).  Composition of an equivalence $\text{Rep}(G)\to\text{Rep}(G')$ with the embedding functor of $\text{Rep}(G')$  is a relaxed tensor functor $\text{Rep}(G)\to \text{Hilb}$.

For the spectral functor ${L}$ of an ergodic action, isometries making it quasitensor    are given 
by $\tilde{{L}}_{u,v}(\sum\psi_i\otimes c_i)\otimes(\phi_j\otimes d_j)=\sum
 (\psi_i\otimes\phi_j)\otimes d_jc_i$. Note that the spectral space of a non spectral representation is trivial, hence, quasitensor functors, unlike the relaxed tensor ones, may take a nonzero object to the zero object.

A quasitensor functor $(\mu,\tilde{\mu})$ preserves conjugates, in the sense that 
if   a nonzero object  $u$ of ${\cal A}$ has a conjugate defined by arrows $R\in(\iota,\overline{u}\otimes u)$ and $\overline{R}\in(\iota,{u}\otimes
  \overline{u})$   then  $\mu_{\overline{u}}$ is a conjugate of
 $\mu_u$, defined 
 by $\hat{R}=\tilde{\mu}_{\overline{u}, u}^*\circ \mu(R)$, $\hat{\overline{R}}=\tilde{\mu}_{{u}, \overline{u}}^*\circ \mu(\overline{R})$.
  This   is a straightforward consequence of the axioms.  
 The property of conservation of conjugates of a quasitensor functor $\mu$ implies   $d(\mu_u)\leq\|\hat{R}\|\|\hat{\overline{R}}\|\leq d(u)$.

In particular, for the spectral functor of an ergodic action we may associate
 an antilinear invertible map $J_u:L_u\to L_{\overline{u}}$ to a solution $(R_u,\overline{R}_u)$ of the conjugate equations of a   representation $u$ of $G$ by
  $$\sum_k J_uT_k\otimes T_k=\tilde{L}_{\overline{u}, u}^*\circ L(R_u),$$ where $T_k$ is an orthonormal basis of $L_u$. 
 The above relation between intrinsic dimensions  becomes a  multiplicity bound of an irreducible representation in the spectum of the related ergodic action. Indeed, the scalar $\|\hat{R}\|\|\hat{\overline{R}}\|$  arising   from the  solution  of the conjugate equations of an irreducible representation $u$, reduces to the 
 quantum multiplicity $q-\text{mult}(u)$ of \cite{BDV}, while the intrinsic dimension of $L_u$ is just the integral dimension of $L_u$, i.e. the 
 ordinary multiplicity $\text{mult}(u)$ of $u$ in the action. The above estimate says 
 $\text{mult}(u)\leq\text{q-mult}(u)\leq d(u)$
 a fact shown in \cite{BDV} refining the inequality $\text{mult}(u)\leq d(u)$ previously obtained by  \cite{Boca}.
 In the classical case, i.e. when $G$ is a compact group, the quantum dimension of a representation is just its integral dimension, and we thus in turn recover the classical result of H\o egh-Krohn Landstad and Stormer that $\text{mult}(u)\leq \text{dim}(H_u)$.
 
 We shall need the following facts proved in  \cite{PR}.

 \medskip
 
 \noindent{\bf 2.2.  Lemma.} {\it Let $G$ be a compact quantum group acting ergodically on a unital $C^*$--algebra with spectral functor $(L,\tilde{L})$. If $J_u$ and $J_v$ are associated to solutions $j_u$ and $j_v$ of the conjugate equations for representations  $u$ and $v$ of $G$ respectively, then 
 \begin{description}
 \item{\rm a)} $J_{u\otimes v} \tilde{L}_{u, v}=\tilde{L}_{\overline{v}, \overline{u}}J_v\otimes J_u\Theta$, where $J_{u\otimes v}$ is associated to the tensor product solution    for $u\otimes v$ and $\Theta: L_u\otimes L_v\to L_v\otimes L_u$ is the flip map,
  \item{\rm b)}
 for any $A\in(u,v)$, $L(j_vAj_u^{-1}) J_u=J_vL(A)$.
 \end{description}
 }\medskip

Recall that conversely, given a  
 quasitensor functor $(\mu,\tilde\mu):\text{Rep}(G)\to\text{Hilb}$, we may associate a maximal ergodic 
 action of $G$ in the following way. 
 Form the linear space ${}^\circ{\cal C}_\mu:=\sum_{u\in\text{Rep}(G)} \overline{\mu_u}\otimes_{\text{Rep}(G)} H_u$, where $\otimes_{\text{Rep}(G)}$ indicates a suitable tensor product treating  the arrows of ${\text{Rep}(G)}$ as scalars. The algebraic operations are defined by, dropping the indices and the natural transformation,
  $$(\overline{k}\otimes\psi)(\overline{k'}\otimes\psi')=(k\otimes k')\otimes(\psi\otimes\psi'), \quad\quad  (\overline{k}\otimes\psi)^*=\overline{Jk}\otimes j^{-1*}\psi.\eqno(2.5)$$ The action of the quantum group on  each subspace $\overline{\mu_u}\otimes H_u$ is the tensor product of the trivial action on the  first factor  and the representation $u$ on the second. This action is ergodic.
 
   Most importantly, the   linear functional $\omega$ on ${}^\circ{\cal C}_\mu$ which annihilates each subspace
 $\overline{\mu_u}\otimes H_u$, $u\in\hat{G}$, $u\neq\iota$ and takes $I$ to $1$ is a positive and  faithful state, it is the unique state invariant under the action  (i.e. $\omega\otimes\text{id}\circ\alpha=\omega$) . Therefore ${}^\circ{\cal C}_\mu$ has a $C^*$--norm. It turns out that the maximal $C^*$--norm is finite. We thus have at 
 our disposal two possible  completions of ${}^\circ{\cal C}_\mu$, which are different in general,  the  completion in the maximal $C^*$--norm, denoted ${\cal C}_\mu$, and 
the completion   in  the norm provided by the GNS representation $\pi_\omega$ of the invariant state, called the  reduced completion.
 
The $G$--action clearly extends to the maximal completion. In the case of the reduced completion, note that the action of $G$ on ${}^\circ{\cal C}_\mu$ is in fact only an action of the dense spectral Hopf $^*$--subalgebra. This action extends to
an action of the reduced compact quantum group $G_{\text{red}}$ and in turn it lifts to a normal action of the Hopf--von Neumann algebra generated by  the regular representation of $G$ on the von Neumann algebra $\pi_\omega({\cal C})''$. In all these cases, the extended action is
  ergodic, see  Theorem 2.5 of \cite{Wang}.
\medskip

 \section{ Modular theory of the invariant state}  In view of the duality theorem recalled in the previous section, we may and shall think of the spectral functor of an ergodic $C^*$--action of a compact quantum group $G$ as an abstract quasitensor functor $(\mu,\tilde{\mu}):\text{Rep}(G)\to\text{Hilb}$. Correspondingly,
 we shall represent the dense spectral subalgebra with generators and relations given in $(2.5)$.
 
If $a=\overline{k}\otimes\psi$, $b=\overline{k'}\otimes\psi'$ have support in the irreducible representations $u$ and $v$ respectively then 
$$\omega(a^*b)=\delta_{u,{v}}\|R_u\|^{-2} (k', {J_u}^*J_uk)  (\psi,\psi'),\eqno (3.1)$$
$$\omega(ba^*)=\delta_{u,v}\|\overline{R}_u\|^{-2}(k',k)(\psi,({j_u}^* j_u)^{-1}\psi'),\eqno(3.2)$$
see Sect. 8 in \cite{PR} for explicit computations. These formulas may be used to derive  modular properties of the invariant state, in turn generalizing the corresponding properties of the Haar state \cite{WoronowiczCMP, Boca}.  
In fact, for every irreducible representation $u$ choose a standard solution
$(R_u,\overline{R}_u)$.   We have a densely defined multiplicative map such that on $\overline{\mu_u}\otimes H_u$,
$$\sigma_{-i}(\overline{k}\otimes\psi)=\overline{(J_u^*J_u)k}\otimes j_u^*j_u\psi$$ 
with inverse
$$\sigma_i(\overline{k}\otimes\psi)=\overline{(J_u^*J_u)^{-1}k}\otimes (j_u^*j_u)^{-1}\psi.\eqno(3.3)$$ 
Since a different standard solution is of the form $Uj_u$, with $U$ unitary, the associated 
$J_u$ changes into $\mu(U)J_u$, hence $\sigma_{-i}$ and its inverse do not change.
A quick computation shows that the KMS property holds,
$$\omega(\sigma_{-i}(b)a^*)=\omega(a^*b).$$
We may choose $j_{\overline{u}}=j_{u}^{-1}$, which implies $J_{\overline{u}}=J_{u}^{-1}$. It follows that
$\sigma_{-i}(a^*)=\sigma_i(a)^*$. We collect the conclusions of the above discussion, a slight refinement of a  result of \cite{WoronowiczCMP, Boca} .
\medskip

\noindent{\bf 3.1. Theorem.} {\it The invariant state of an ergodic action of a compact quantum group $G$ on a unital $C^*$--algebra satisfies the KMS condition on the dense spectral subalgebra. It is a trace if and only if
 for any spectral irreducible representation $u$ of $G$,
\begin{description}
\item{\rm a)} $d(u)=\text{\rm dim}(u)$,
i.e. the antilinear $j_u$ defining a standard solution is antiunitary, and
\item{\rm b)} the associated $J_u$ is  antiunitary  (hence $q-\text{\rm mult}(u)=\text{\rm mult}(u)$).
\end{description}}\medskip
\medskip

\noindent Recall that the condition $d(u)=\text{dim}(u)$ for all $u$ means precisely that the quantum group is of Kac type, i.e. it has involutive coinverse. 
 For example, quantum dimensions of 
the real deformations $G_q$, $0<q<1$ of classical compact Lie groups are known to be strictly larger than the corresponding integral dimensions, hence the invariant state
of any ergodic $C^*$--algebra under the action of any of these quantum groups is not a trace.
\medskip

\noindent{\it Almost periodicity of  the invariant state.}
If $({\cal C},\alpha,G)$ is an ergodic action of a compact quantum group, the dense spectral subalgebra of ${\cal C}$ is the domain
of a    one parameter group of $^*$--automorphisms, the modular group, given by $$\sigma_t(\overline{k}\otimes\psi)=
\overline{({J_u}^*J_u)^{it}k}\otimes (j_u^*j_u)^{-it}\psi.$$
 This group extends 
to a one parameter automorphism group of the maximal completion ${\cal C}_\mu$. Moreover, it extends to $\pi_\omega({\cal C}_\mu)$ or $M:=\pi_\omega({\cal C}_\mu)''$ as well since it leaves $\omega$ invariant. In this subsection we only consider the extension to the von Neumann algebra $M$.
 Note that
 the cyclic vector $\Omega$ associated to the GNS representation of   
 $\omega$ is separating for   $M$. This may be shown with arguments similar to those of Sect. 4 of \cite{DLRZ}. Hence, by the KMS property, $\sigma_i$ becomes the restriction of the modular operator $\Delta_\omega$ associated to $\omega$
 under the canonical inclusion ${}^\circ{\cal C}_\mu\to L^2({\cal C},\omega)$.
 Recall that Connes defined a normal faithful state $\phi$ on a von Neumann algebra
 to be {\it almost periodic} if   $\Delta_\phi$ is diagonal.
In the case of an ergodic action,  the expressions for the inner product $(3.1)$ and  for  $\sigma_i$, $(3.3)$,  show that the cyclic  state of $M$ thus obtained is always almost periodic.
 Moreover, the point spectrum of $\Delta_\omega$ is completely determined by 
 the eigenvalues if 
$j_u^*j_u$ and $J_u^*J_u$ for $u$ describing a complete set of irreducible spectral representations. While the spectrum of $j_u^*j_u$ is 
 a structural property
of the quantum group, that  of $J_u^*J_u$,   depends on the ergodic action, and this is 
the most mysterious part. Indeed, although $J_u$ is explicitly associated with $j_u$, 
we can not infer that properties of $j_u$ pass to $J_u$. For example, we may have
$j_u$ antiunitary for all spectral $u$ but $J_u$ not antiunitary, or, in other words, we may have ergodic actions of compact quantum groups of Kac type on a unital $C^*$--algebra ${\cal C}$ for which $M=\pi_\omega({\cal C})''$ is a type III factor \cite{Wang}. Similarly, we may have $j_u$ not antiunitary for all spectral $u$ but $J_u$ always antiunitary, as in the examples arising from subfactors, described in detail in \cite{PRsubfactors}.
\medskip

We conclude this section with a few more remarks on the modular theory of  $({\cal C}, \alpha, G)$.   Let $\text{Sp}(\Delta_\omega)$ ($\text{Sp}_p(\Delta_\omega)$) denote the 
  spectrum (point spectrum) of   $\Delta_\omega$. The following fact may be known,
  a proof is included for convenience.
 
  \medskip
  
  \noindent{\bf 3.2. Theorem.} {\it If an ergodic $C^*$--action of a compact quantum group $G$  on ${\cal C}$ admits a spectral irreducible representation $u$ of $G$ such that $d(u)>\text{\rm dim}(u)$   then  $\text{\rm Sp}_p(\Delta_\omega)\neq \{1\}$.
 If in addition both $M=\pi_\omega({\cal C})''$ and  the centralizer $M_\omega$ are  factors, then $M$ is of type {\rm III}${}_\lambda$ with $0<\lambda\leq1$.  }\medskip
  
  \noindent{\it Proof}  
  If $u$ is a
   spectral irreducible representation for which
  $j_u$ is not antiunitary,  then $j_u^*j_u$ has an eigenvalue $<1$ and another $>1$ since $j_u$ is standard. Hence   $\Delta_\omega$ has an eigenvalue $\neq1$ on $\overline{\mu_u}\otimes H_u$. If $M$ and $M_\omega$ are factors, it is well known that
  $S(M)=\text{Sp}(\Delta_\omega)$, with $S(M)$ the Connes invariant (Cor. 3.2.7 in
  \cite{ConnesIII}).
   Hence $S(M)\neq\{1\}$ and $S(M)\neq\{0,1\}$, so $M$ is not semifinite or of type {\rm III}${}_0$.
   \medskip  

In the next section we shall give  general lower and upper bounds for the eigenvalues of the modular operator $\Delta_\omega$, depending only on the quantum group. We shall derive a general lower bound for the possible parameters $0<\lambda<1$ such that  $M=\pi_\omega({\cal C})''$ is a factor of type III${}_\lambda$  with a factorial centralizer $M_\omega$.

\medskip

\section{An analogue of the finiteness theorem for CQG}

We have thus seen that the point spectrum of the modular operator $\Delta_\omega$ associated to an ergodic action is completely determined by the spectra of the positive matrices
$(J_u^*J_u)\otimes (j_u^*j_u)$ associated to the spectral irreducible representations $u$ of the quantum group.

In this section we derive a general estimate for the eigenvalues of the operators $J_u^*J_u$. Our estimate involves the growth rate of the quantum dimension of $u$. It  allows to reproduce  the finiteness result of \cite{HLS}, and also   to derive     other modular properties of the invariant state.

We define the growth rate  of the intrinsic dimension of an object  $u$ of a tensor $C^*$--category,
$$D_{u,n}:=\max{\{d(v), \text{ irreducible subobjects $v$ of } u^{\otimes n}\}},$$
 $$D_u:=\lim_n (D_{u,n})^{1/n}.$$
This limit always exists since $D_{u,n}$ is a submultiplicative sequence.  Note that $D_u\leq d(u)$, $D_u=D_v$ if $u$ and $v$ are  equivalent, and $D_{\overline{u}}=D_u$.
We say that  the intrinsic dimension of  $u$   has subexponential  (exponential) growth if $D_u=1$ ($D_u>1$).
Note that if $w$ is a subrepresentation of some tensor power $u^{\otimes k}$ of $u$,
$D_w\leq D_u^k$. Hence, if the tensor powers of $u$ contain all the irreducibles and if $d(u)$ has subexponential growth, so does every irreducible.\medskip

\noindent{\bf 4.1. Example.} Let  $G$ be a compact group. In this case the quantum dimension of any   representation $u$  is just the integral dimension of the corresponding Hilbert space. It is known that this dimension always has polynomial, and hence subexponential, growth.  This fact  relies on Weyl's dimension formula and played an important role in 
the original proof of the finiteness theorem of \cite{HLS}.  
\medskip

 Let $({\cal C},\alpha)$ be an ergodic action of a compact quantum group $G$ on a unital $C^*$--algebra.
We shall also need to consider the growth rate of integral multiplicities of irreducible representations of $G$. Given a representation $u$, set
$$\text{Mult}_{u,n}:=\max\{\text{mult}(v), \text{irreducible subrepresentations $v$ of $u^{\otimes n}$}\},$$
$$\text{Mult}_u:=\liminf_n \text{Mult}_{u,n}^{1/n},$$
 which enjoys properties similar to those of $D_u$. In the particular case of the translation action of $G$ over itself, integral multiplicities $\text{mult}(u)$ reduce to  integral dimensions $\text{dim}(u)$. We shall accordingly denote the corresponding growth rate by $\text{Dim}_u$ and refer to it as the growth rate of integral multiplicities. In both cases,
 we have the notion of subexponential or exponential growth. A comparison between the various growth rates introduced
 may be easily derived,
 $$ \text{Mult}_u\leq D_u,\quad \text{Dim}_u\leq D_u.$$
 \medskip

\noindent{\bf 4.2. Example.}
Let $G$ be a classical compact Lie group and $G_q$ the associated deformed
compact matrix quantum group by a positive parameter $0<q<1$. The integral dimensions of irreducibles have subexponential growth as they are the same as in the classical case. On the other hand, the quantum dimensions of irreducibles have exponential growth rate. This last assertion is known. However, it may also be derived from the following proposition. Indeed, $D_u=1$ would imply $\lambda_u=\Lambda_u=1$, hence $j_u$ antiunitary.
\medskip

\noindent{\bf 4.3. Proposition.} {\it For any representation $u$ of a compact quantum group, $$D_u^{-1}\leq\lambda_u\leq\Lambda_u\leq D_u.$$
If $\text{\rm Dim}_u=1$ then
the first and last inequalities are equalities.}\medskip

\noindent{\it Proof} The middle inequality being obvious, it suffices to show that $\Lambda_u\leq D_u$,
as $\lambda_u\geq D_u^{-1}$ follows   passing to the conjugate representation.
The $n$--th tensor product solution $j_{u^{\otimes n}}=(j_u\otimes\dots\otimes j_u)\theta_n$ for $u^{\otimes n}$, with $\theta_n$ a suitable permutation operator,
is standard if $j_u$ is.
 If $v$ is an irreducible subrepresentation of $u^{\otimes n}$ and $j_v$ is a standard solution for $v$, $d(v)=\text{Trace}(j_v^*j_v)\geq\Lambda_v$.
Hence $$D_{u,n}\geq\max\{\Lambda_v, v \text{ irreducible subrepresentation of } u^{\otimes n}\}=$$
$$\Lambda_{u^{\otimes n}}=\Lambda_u^n.$$ 
On the other hand $d(v)\leq\Lambda_v\text{dim}(v)$ implies
$$D_{u,n}\leq\Lambda_{u^{\otimes n}}\text{Dim}_{u,n}=\Lambda_u^n\text{Dim}_{u,n},$$
hence $D_u\leq \Lambda_u$ if the integral dimension of $u$ has subexponential growth rate. 
\medskip

\noindent{\it Remark}   This proposition may be used to
explain why the spectrum of the associated matrices
$j_u^*j_u$ has a symmetric shape for the deformed compact matrix quantum groups $G_q$. 
\medskip

For the Wang-Van Daele quantum groups $A_o(F)$ and $A_u(F)$, we follow 
the notation of \cite{Banica}. We shall always normalize $F$ so that $\text{Trace}(F^*F)=\text{Trace}((F^*F)^{-1})$. Recall that for $A_o(F)$, the matrix
 $F$ is required to satisfy $F\overline{F}=\pm 1$.
 \medskip

\noindent{\bf 4.4. Example.} If $u$ is the fundamental representation of $A_u(F)$, with $\text{rank}(F)\geq2$, all tensor powers $u^{\otimes n}$ are irreducible \cite{Banica}, showing that $D_{u,n}=d(u)^n$ by multiplicativity of quantum dimension. Hence $D_u=d(u)=\text{Trace}(F^*F)$ is the largest possible value. Similarly, $\text{Dim}_u=\text{dim}(u)$. In particular,
the extreme inequalities in Prop. 4.3 are always strict 
(the first examples being 
  $A_u(m)$, $\lambda_u=\Lambda_u=1$ but $D_u=\text{Dim}_u=m$),
while the middle inequality is generically strict for $\text{rank}(F)\geq3$, but it is an equality for $\text{rank}(F)=2$ due to $\text{Trace}(F^*F)=\text{Trace}((F^*F)^{-1})$. 
\medskip

\noindent{\bf 4.5. Example.} If $u_1$ is the fundamental representation  of  $G=S_qU(2)$, for a nonzero $|q|\leq 1$,
$|q|^{1/2}R=\psi_1\otimes\psi_2-q\psi_2\otimes\psi_1$ is a standard solution, hence $j_{u_1}^*j_{u_1}=\text{diag}(|q|, |q|^{-1})$. Let $u_r$ be the unique irreducible $r+1$--dimensional representation. The Clebsch--Gordan rule $u_1\otimes u_r\simeq u_{r-1}\oplus u_{r+1}$ gives 
$j_{u_r}^*j_{u_r}=\text{diag}(|q|^r, |q|^{r-2},\dots,|q|^{-r})$. Hence $\lambda_{u_r}=\Lambda_{u_r}^{-1}=D_{u_r}^{-1}=|q|^r$.

\medskip

\noindent{\bf 4.6. Example.} Consider the quantum group $A_o(F)$.  It is well known that the representation categories of $A_o(F)$ and $S_{\mp q}U(2)$ are tensor equivalent if $q>0$ is defined by $q+\frac{1}{q}=\text{Trace}(F^*F)$. Hence 
   $D^{A_o(F)}_{u_r}=q^{-r}$ by 4.5. Note that $\lambda_{u_r}=\Lambda_{u_r}^{-1}$.
   This follows from the validity for $r=1$ and the Clebsch-Gordan rule.
A computation shows that for every $F$ and   for the fundamental representation $u_1$, the extreme inequalities
in Prop. 4.3 are strict iff $\text{rank}(F)\geq3$.

\medskip

\noindent{\bf 4.7. Example.}
For $\text{rank}(F)\geq3$,    the integral dimension  of the fundamental representation $u_1$  of  $A_o(F)$ has exponential growth.  This may be seen in the following way. Independently of the matrix $F$, irreducible representations of $A_o(F)$ 
satisfy the same fusion rules as those of $SU(2)$ \cite{Banica}. Hence, denoting with the same symbol the corresponding irreducible representations,   the integral dimensions
are determined by the Clebsch--Gordan rule, $\text{dim}(u_{r+1})=\text{dim}(u_1)\text{dim}(u_r)-\text{dim}(u_{r-1})\geq 2\text{dim}(u_r)$, hence
$\text{Dim}_{u_1, r}=\text{dim}(u_r)\geq 2^{r-1}\text{dim}(u_1)$.
\medskip

The following is our main result. \medskip

\noindent{\bf 4.8. Theorem.} {\it For any spectral irreducible representation $u$ of an ergodic action of a compact quantum group on a unital $C^*$--algebra,
\begin{description}
 \item{\rm a)}
 $D_u^{-1}\leq J_u^*J_u\leq D_u,$
 where $J_u$ is associated to a standard solution $j_u$ for $u$.
 \item{\rm b)}
  If the spectral functor is relaxed tensor and if $\text{\rm Mult}_u=1$  then $\frac{1}{D_u}$ and $D_u$ are respectively the smallest and the largest eigenvalues of $J_u^*J_u$. 
 \end{description}
 }\medskip

\noindent{\it Proof} a) The first inequality follows from the second
applied to $\overline{u}$.
By Lemma 2.2 a), for any positive integer $n$, if $J_u$ is associated to any solution $j_u$ for $u$,
$$\|J_u^*J_u\|^n=\|J_u^*J_u\otimes\dots\otimes J _u^*J_u\|=$$
$$\|\tilde{\mu}_{u,\dots,u}^*J_{u^{\otimes n}}^*J_{u^{\otimes n}}\tilde{\mu}_{u,\dots,u}\|\leq\|J_{u^{\otimes n}}^*J_{u^{\otimes n}}\|,$$  
where $J_{u^{\otimes n}}$ is associated to the $n$--th tensor product solution $j_{u^{\otimes n}}=(j_u\otimes\dots\otimes j_u)\theta_n$ for $u^{\otimes n}$.
 Consider a complete reduction of $u^{\otimes n}$ into a direct sum of irreducible 
representations $v$ and let $S_v\in(v, u^{\otimes n})$ be the isometry associated to $v$.
We may compute $\|J_{u^{\otimes n}}^*J_{u^{\otimes n}}\|$   as the norm of the positive operator--valued matrix 
$(\mu(S_v^*)J_{u^{\otimes n}}^*J_{u^{\otimes n}}\mu(S_w))_{v,w}$.
 By lemma 2.2 b),
$\mu(j_{u^{\otimes n}}{S_v}j_{v}^{-1})J_v=J_{u^{\otimes n}}\mu(S_v)$, with respect to any solution $j_v$ for $v$, with associated $J_v$. 
Now we fix  standard solutions for $u$ and all the $v$. Since $j_{u^{\otimes n}}$ is standard,  we may find $\{S_v\}$ such that $j_{u^{\otimes n}}S_vj_{v}^{-1}=:\overline{S}_v$ are pairwise orthogonal
isometries.
Hence
$$\mu(S_v^*)J_{u^{\otimes n}}^*J_{u^{\otimes n}}\mu(S_w)=J_{v}^*\mu(\overline{S}_v^*\overline{S}_w)J_w=\delta_{v,w}J_v^*J_w.$$ 
Combination with the previous estimate gives,  
$$\|J_u^*J_u\|\leq \|\text{diag}_v(J_v^*J_v)\|^{1/n}=(\max_v\{\|J_v^*J_v\|\})^{1/n}\leq$$
$$
(\max_v\{\text{Trace}(J_v^*J_v)\})^{1/n}=(\max_v\{\|\hat{R_v}\|^2\})^{1/n}\leq
(\max_v\{\|{R_v}\|^2\})^{1/n}=
D_{u,n}^{1/n}.$$

b) If the spectral functor is relaxed tensor, all the $\tilde{\mu}$ are unitary, hence the first and the last inequalities are equalities. On the other hand the estimate $\max_v\{\text{Trace}(J_v^*J_v)\}\leq \max_v\{\|J_v^*J_v\|\}\text{Mult}_{u,n}$ shows that $\|J_u^*J_u\|=D_u$ if the integral multiplicity of $u$ has subexponential growth. 
\medskip

\noindent{\bf 4.9. Corollary.}  \cite{HLS} {\it  The invariant state of an ergodic action of a compact group on a unital $C^*$--algebra is a trace.}\medskip

\noindent{\it Proof}  As recalled above, for all $u$, $D_{u,n}$ has polynomial growth. Hence $J_u$ are all antiunitary by Theorem 4.8. On the other hand the same holds in the classical case for standard solutions $j_u$ of $G$--respresentations,   hence we may apply Theorem 3.1.
\medskip

We next discuss examples satisfying b) with $D_u>1$ and which are not translation actions.\medskip

\noindent{\bf 4.10. Example.} Let $G_q$ be as above a deformed classical compact Lie group, and let  $\Phi:\text{Rep}(G')\to\text{Rep}(G_q)$ be a tensor equivalence, with $G'$  another compact matrix quantum group (e.g.  $G=SU(2)$ and $G'=A_o(F)$).
The composition $\mu: \text{Rep}(G')\to\text{Rep}(G_q)\to\text{Hilb}$ with the embedding functor of $\text{Rep}(G_q)$ is a relaxed tensor functor. 
Consider the maximal ergodic action associated to $\mu$. The integral multiplicity of an irreducible representation $u$ of $G'$ is the integral dimension of $\Phi(u)$, which is the same as in the classical case,
hence $\text{Mult}_u^{G'}=\text{Dim}_{\Phi(u)}^{G_q}=1$. 
 Note that these are just  the inverses of the tensor equivalences of \cite{BDV}.
 \medskip

\noindent{\bf 4.11. Corollary.} {\it If a  compact quantum group $G$ of Kac type admits an ergodic action on a unital $C^*$--algebra with non-tracial invariant state then  the integral dimension of some spectral irreducible representation  of $G$ has exponential growth.}\medskip

\noindent{\it Proof} Since $G$ has involutive coinverse, every standard solution $j_u$
is antiunitary. If all  the quantum dimensions of spectral irreducibles $u$ had subexponential growth, the associated $J_u$ would be antiunitary by Theorem 4.8.,
hence the  invariant state would be a trace by Theorem 3.1.
\medskip

\noindent{\bf 4.12. Corollary.} {\it $S_{-1}U(d)$ acts ergodically only on tracial $C^*$--algebras.}\medskip

\noindent{\bf 4.13. Corollary.} {\it Let $G$ be a compact quantum group acting ergodically on ${\cal C}$ and admitting a spectral irreducible representation $u$ s.t. $d(u)>\text{\rm dim}(u)$. Assume that 
   $M:=\pi_\omega({\cal C})''$ and $M_\omega$ are factors.
   \begin{description}
   \item{\rm a)} If $M$ 
   is of type
 {\rm III}${}_\lambda$,  $0<\lambda<1$, then 
$$\lambda\geq\sup\{ \frac{\min\{\lambda_u, \Lambda_u^{-1}\}}{D_u}, u \text{ spectral irr. s.t. } d(u)>\text{\rm dim}(u)\},$$
\item{\rm b)} if the above supremum is $1$, $M$ is of type {\rm III}${}_1$.
\end{description}

}\medskip

\noindent{\it Proof} a) Let $u$ be a spectral irreducible representation such that $d(u)>\text{dim}(u)$. As argued in Theorem 3.2., the modular group of $\omega$ does not act trivially on the spectral subspace $\overline{\mu_u}\otimes H_u$. On the other hand, the smallest and largest eigenvalues of the restriction of the modular 
operator $\Delta_\omega$ to $\overline{\mu_u}\otimes H_u$ are bounded below  by $\frac{1}{D_u\Lambda_u}$ and above by
$\frac{D_u}{\lambda_u}$. Taking into consideration the fact that the eigenvalues of $\Delta_\omega$ belong to $S(M)$, we see that either
$\frac{1}{D_u\Lambda_u}\leq\lambda$ or $\lambda^{-1}\leq \frac{D_u}{\lambda_u}$.
b) If the supremum is 1, by a), $M$ is not of type III${}_\lambda$ for any $0<\lambda<1$.
Hence $M$ must be of type III${}_1$ by Theorem 3.2.
\medskip

\noindent{\it Remark}  If  the quantum group is of Kac type, under the same assumptions of the previous corollary, we similarly derive   $\lambda\geq\frac{1}{n}$, with $n$ the minimal dimension of a spectral irreducible representation for which the restriction of the modular group on the associated spectral subspace is non-trivial. 
\medskip

\noindent{\bf 4.14. Example.} The lower bound in the above remark is realized in
Wang's example of  ergodic action of $A_u(n)$ on the type III${}_{\frac{1}{n}}$ 
factors \cite{Wang}. More precisely, in  this case the factor is the von Neumann completion of the Cuntz algebra
${\cal O}_n$ in the GNS representation of the canonical state and the modular group is 
the group describing  the ${\mathbb Z}$--gradation. This automorphism group acts non-trivially on the generating Hilbert space of isometries. This Hilbert space carries the fundamental representation of the quantum group. 
This example shows that $D_u$, together with triviality of the modular theory of $A_u(n)$, explain completely the factorial type.\medskip

\noindent{\it Remark} If $\lambda_u=D_u^{-1}={\Lambda_u}^{-1}$,
 e.g. under the condition of Prop. 4.3,  the lower bound for $\lambda$ becomes simply $\sup\{\lambda_u^2, \ u: d(u)>\text{dim}(u)\}$. 
In particular, for $S_qU(2)$, with $0<|q|<1$, under the same assumptions of the previous theorem, the possible parameters $\lambda$  satisfy $\lambda\geq |q|^{2r}$, where
$r$ is the smallest strictly positive integer such that $u_r$ is spectral. Hence, for a given $q$,  small values of $\lambda$ are possible only if the spectrum of the action has gaps.
\medskip

The following example shows  that the assumption of subexponential  growth of integral multiplicities can not be removed in Theorem 4.8 b).\medskip 

\noindent{\bf 4.15. Example.}
For $A_o(F)$, with $F$ of rank $\geq 3$,  taking into account 
Example 4.6 and $\lambda_{u_r}=\Lambda_{u_r}^{-1}=\|F\|^{-2r}$,  we see that $\lambda\geq(\frac{q}{\|F\|^2})^r$, where $q$ is defined by  $\text{Trace}(F^*F)=\text{Trace}((F^*F)^{-1})=q+\frac{1}{q}$ and  $u_r$  is again the first spectral irreducible.
In particular,
 consider   the translation action of $A_o(F)$ over itself.
 Vaes and Vergnioux  have shown, among other things, that if $\text{Trace}(F^*F)\geq\sqrt{5}\|F\|^2$ then the Haar state $h$ is factorial  and the associated von Neumann algebra is a full factor. 
 Moreover, if the spectrum of $(F^*F)^{-1}\otimes F^*F$ generates the subgroup $\{\lambda^n, n\in{\mathbb Z}\}$, for some  $0<\lambda<1$, $\pi_h(A_o(F))''$ is a factor of type III$_{\lambda}$  \cite{VV}. 
 The fundamental representation is spectral.
 In this example our lower bound 
 $\lambda\geq\frac{q}{\|F\|^2}$ is not optimal  already for 
$F^*F=\text{diag}(\lambda, 1, \lambda^{-1})$. 
\medskip

\noindent{\it Acknowledgements.} I would like to thank Alessandro Fig\`a-Talamanca for a conversation related to this note.

\noindent Claudia Pinzari, 
Dipartimento di Matematica, Sapienza Universit\`a di Roma,
00185--Roma, Italy; e-mail: pinzari@mat.uniroma1.it

\end{document}